\documentclass{amsart}

\let\ep=\varepsilon
\def\hess{\mathop{i\partial {\overline\partial}}\nolimits }
\def\card{\mathop{card}}

\def\n{|\!|}

\def\om{\omega}
\def\cb{\mathbb{C}}
\def\nb{\mathbb{N}}

\def\rb{\mathbb{R}}
\def\bb{\mathbb{B}}

\def\fc{\mathcal{F}}

\newtheorem{theorem}{Theorem}
 \newtheorem{prop}{Proposition}
\newtheorem{definition}{Definition}

\begin{document}
\title{Interpolation in non-positively curved K\"ahler manifolds}
\author{Christophe Mourougane}

\maketitle

{We   extend  to   any  simply   connected  K\"ahler   manifold  with
non-positive  sectional   curvature  some   conditions  for
interpolation  in   $\mathbb{C}$  and  in  the  unit   disk  given  by
Berndtsson, Ortega-Cerd\`a  and Seip.  The main tool  is a comparison
theorem  for the  Hessian in  K\"ahler geometry  due to  Greene, Wu and
Siu, Yau.}

\section{Introduction}
\label{sec:intro}
Let $(X,\om)$ be a K{\"a}hler manifold, $(L,h)\to X$ be
a holomorphic hermitian line bundle on $X$, and  $\Lambda$ a discrete 
set of points in $X$. 
Consider
$$\fc _{\om,h} :=\{f\in H^0(X,L) / \int_{X} \n f\n_h^2 dV_\om
<+\infty\}$$ and  
$$l^2 _{\Lambda, h} :=\{a\in\Gamma (\Lambda, L) / \sum_{p\in\Lambda}
\n a(p)\n_h^2 <+\infty\}$$

We say that the discrete set $\Lambda$ is interpolating for $\fc_{\om,h}$, 
if for every  $a$ in $l^2 _{\Lambda, h}$ there is a section 
$f$ in $\fc _{\om,h}$ such that $f _{\vert \Lambda} =a$.

\bigskip

In the complex plane with its usual flat metric,
there is a complete characterisation of
interpolating sets for pluri-sub-harmonic weights comparable to $|z|^2$.

\textbf{Theorem}[Berndtsson and Ortega Cerd{\`a}~\cite{bo}, Ortega-Cerd{\`a}   
and Seip~\cite{os}]

Consider the trivial line bundle on $\cb$
with metric $h_\Phi:=\vert ~ ~ \vert^2 e ^ {-\Phi}$ where
 $\Phi$ is a pluri-sub-harmonic function on $\cb$ whose Laplacian is 
uniformly bounded. 
A set $\Lambda$ in $\cb$ is interpolating for $\fc _{h_\Phi}$ 
if and only if it is uniformly separated and 
$$\exists \rho>0, \ep >0\ \  \forall z\in \cb,  \ \ \Delta \Phi (z)\geq 
\frac{\card B(z,\rho)\cap \Lambda}{\rho^2} +\ep .$$
Here $B(z,\rho)$ is the Euclidean ball of centre $z$ and radius $\rho$.

\bigskip A simply-connected complete Riemannian manifold   with
non-positive sectional curvature is called a 
\textbf{Cartan-Hadamard} manifold. An hermitian complex manifold is
called Cartan-Hadamard if its underlying Riemannian space is
Cartan-Hadamard. 

We prove 
\begin{theorem}\label{theo1}
Let $(X,\om)$ be a K{\"a}hler Cartan-Hadamard manifold with sectional
curvature bounded from below by $-k^2$. 
Let $\Lambda $ be a discrete set of
points in $X$ and $(L,h)\to X$ be an holomorphic hermitian
line bundle on $X$ almost analytic 
around $\Lambda$  (see definition~\ref{defi:almost analytic}).

 If $\Lambda$ is a uniformly $\om$-separated set such that
$\exists \rho >0, \ep >0,   \forall z\in X$ ,  
$$ic_h(L)+ ricci(\om ) \geq  
 n\frac{\card B_\om (z,\rho)\cap \Lambda}{\rho^2}\left( 1+k\rho \coth
   k\rho\right)\om +\ep\om $$
then it is interpolating for $\fc _{\om,h}$. 
\end{theorem}

 On the other hand, Lindholm~\cite{li} obtained on $\cb ^n$
necessary conditions for interpolation in terms of the Monge-Amp{\`e}re
measure $(\hess\Phi)^n$ of the weight. For $n>1$, there is therefore a
 gap between necessary and sufficient conditions expressed in terms 
of eigenvalues of the curvature.  

\bigskip
Seip also obtained a complete description of interpolating sets in the
unit disk with its hyperbolic measure for special weights $-A\log
(1-\vert z\vert ^2)$ in the trivial line bundle in terms of the
hyperbolic lower density for $\Lambda$.
In the spirit of~\cite{seip},  we define for
 $(X,\om)$  a K{\"a}hler Cartan-Hadamard manifold  with sectional curvature 
bounded from above by $-1/\kappa ^ 2$ , $\Lambda$ a
 discrete set of 
points in $X$, and $x\in X$
$$D_{\Lambda ,\kappa} (x)
:=\sum_{p\in \Lambda \atop d_\om (x,p)\geq 1} -\log\tanh ^2
\left( \frac{d_\om (x,p)}{2\kappa}\right).$$

We prove 
\begin{theorem}\label{theo2} 
Let $(X,\om)$ be a K{\"a}hler Cartan-Hadamard manifold with sectional
curvature bounded from above by $-1/\kappa ^ 2$ and from below by
$-k^2$. 
Assume $\Lambda$ is uniformly $\om$-separated,
and $(L,h)$ is almost analytic around $\Lambda$. Then if  
$$\sup _X  D_{\lambda ,\kappa}
<+\infty \text{ and } ic_h(L)+ ricci (\om ) \geq \ep\om $$
for some positive $\ep$, the set $\Lambda$ is interpolating for $\fc _{\om,h}$.
\end{theorem}

In the next section, we recall some comparison theorem in Riemannian
geometry from ~\cite{g-w} and ~\cite{s-y}. In
section~\ref{sec:separation}, we 
prove that under some assumption relating the complex structure and
the metric structure uniformly on $\Lambda$, the uniform separation
is necessary. The section~\ref{sec:non-positive} is devoted to prove 
theorem~\ref{theo1} and the section~\ref{sec:negative} to prove 
theorem~\ref{theo2}.

Note that using the extension theorem of
Ohsawa-Takegoshi-Manivel~\cite{ma}, one can get interpolation results on
Stein manifolds, expressed in terms of an equation for the set
$\Lambda$. 

In the whole work, we will denote by $C$ any constant which only
depends on quantities fixed before.

\section{Comparison theorems in Riemannian geometry} 
\label{sec:separation}
 In a Cartan-Hadamard manifold the exponential map is at each point a
 diffeomorphism. There is no conjugate points. 
 The distance between two points is achieved by a unique geodesic.
Hence, the distance function from any point is
 smooth. 
 Recall that by a result of Wu (see~\cite{wu} page 100) a 
 K{\"a}hler Cartan-Hadamard manifold is Stein.
 As a consequence of these facts, any holomorphic line bundle on a 
K{\"a}hler Cartan-Hadamard manifold is trivialisable. Nevertheless, 
it is more convenient to state the results, especially the
 admissibility condition, in terms of hermitian line
 bundles instead of pluri-sub-harmonic weights. 

\bigskip
By comparison theorem of Rauch (see for example~\cite{k-n} page 76),
lower bounds on the sectional curvature yields upper bounds on the
Jacobi vector fields and hence on the norm of the differential of the
exponential map.  
By other forms of the comparison theorem (see~\cite{cha} page
118-119), upper bounds on sectional curvature yields lower bounds on
volume. 

We will state the comparison theorem for the Hessian (see~\cite{g-w}
Theorem A page 19) in a form suitable for our purpose. The proof of it 
relies on a precise form of the formula for second order variation
of arc length, and on a relation between real and complex Hessian
on K\"ahler manifolds (see~\cite{wu} page 102).

\textbf{Theorem}
Let $M$ be a K{\"a}hler Cartan-Hadamard manifold of dimension $n$ and 
constant sectional curvature $k$. 
Let $f$ be a non-decreasing function on $[0,+\infty[$. Let $0$ be
a point in $M$. Assume $\hess f(d_M(0,\cdot))=g(d_M(0,\cdot))\om_M$
for a function $g$ on $[0,+\infty[$.

Let $X$ a K{\"a}hler Cartan-Hadamard manifold of dimension $n$ and sectional
curvature bounded from below (resp. from above) by $k$. Let $x_0\in
X$. Then
$$\hess  f(d_X(x_0,\cdot)) \leq g (d_X(x_0,\cdot))\om_X.$$
resp.
$$\hess  f(d_X(x_0,\cdot)) \geq g (d_X(x_0,\cdot))\om_X.$$

\section{On uniform separation}

We generalize the condition on the weight posed in earlier works on
interpolation. 

\begin{definition} Fix $M_2$ in $[0,+\infty [$.
An hermitian line bundle $(L,h)\to U$ on an open set $U$ of $\cb^ n$
 has an $M_2$-regular frame $e$ on $U$ 
 if the weight function $\Phi:= -\log h(e)$ 
can be written as the sum 
$$\sum_{1\leq \alpha\leq N} \vert \sigma_\alpha\vert ^2 + \Phi_{def}$$
where the $\sigma_\alpha$ are holomorphic functions on $U$ 
with oscillation 
$\sup_{x,y\in U} \vert \sigma_\alpha (x)-\sigma_\alpha (y)\vert^2$ bounded
by $M_2$ and 
all second order real derivatives of the deformation weight $\Phi_{def}$
are bounded by $-M_2$ and  $M_2$ on $U$.
\end{definition}

\begin{definition}\label{defi:regular}
An hermitian line bundle $(L,h)\to X$ is said to be regular
around $\Lambda$  if there exists positive numbers
 $r,\lambda, \mu$ and $M_2$
such that for all $p\in\Lambda$, there exists a coordinates chart 
$\Psi_p : (U,0)\subset \cb^n\to (V,p)\subset X$ centered at $p$ such
that 
\begin{itemize}
\item[(i)] $\Psi _p(U)\supset B_\om (p,r)$

\item[(ii)] $(\Psi_p^{\star} L,h)$ has an $M_2$-regular frame $e$ on
  $\Psi_p^{-1}(B_\om(p,r))$. 

\item[(iii)] The differential $d\Psi_p :(TU,\Psi_p^\star \om (p) )\to
  (TV,\om)$ satisfies $\mu \leq \n d\Psi_p \n \leq \lambda$ on
  $\Psi_p^{-1}(B_\om (p,r))$. 
\end{itemize}
\end{definition}

Remark first that \text{(iii)} implies
$$\mu d_{eucl} \leq d_\om\leq \lambda d_{eucl} \ \ \text{ on } B_\om
(p,r)$$
where $d_{eucl}$ is computed with the constant metric $\Psi_p^\star
\om (p) $ on $U$.

In this section, we will prove the following
\begin{theorem}\label{theo:separation} 
Let $X,L$ and $\Lambda$ as in the introduction.
Assume $X$ is a K{\"a}hler Cartan-Hadamard manifold
and that $(L,h)\to X$ is regular around $\Lambda$  .

Then, if $\Lambda$ is interpolating for $\fc _{\om,h}$, it is
uniformly $\om$-separated.
\end{theorem}

\subsection{On regularity}
It can be proved from the fact that $\nb$ is interpolating in $\cb$ 
for the weight $\vert z\vert ^2$ that the non-uniformly separated set
$\{\sqrt{n}, n\in\nb\}$ is interpolating in $\cb$ 
for the weight $\vert z^2\vert ^2$.
This explains our restriction on the metric $h$.

\subsection{Mean value estimates}
According to~\cite{g-w} (Theorem B page 43), for all $p$ in a K{\"a}hler
Cartan-Hadamard manifold, for all holomorphic function $f$ on $B_\om
(p,r)$ and all $x\in B_\om (p,r/2)$,  
$$\vert f(x)\vert^2\leq \frac{1}{vol_{eucl}(B_{eucl}(r/2))}
\int_{B_\om (x,r/2)} \vert f\vert^2 dV_\om .$$

We now assume the hypothesis of theorem~\ref{theo:separation}.
The following estimates relies on the simple inequalities for $x,z$ in
$B_\om (p,r)$

\begin{eqnarray*}
  \vert \sigma_\alpha (x)\vert ^2 +2\Re \left(
\overline{\sigma_\alpha (x)}(\sigma_\alpha (z)-\sigma_\alpha (x))\right)& = &
\vert \sigma_\alpha (z)\vert ^2 
- \vert \sigma_\alpha (z)-\sigma_\alpha (x)\vert ^2 \\
&\geq & \vert \sigma_\alpha (z)\vert ^2 -M_2
\end{eqnarray*}
and 
 with the coordinates $(z_k)$ provided by $\Psi_p$,
\begin{eqnarray*}
  \Phi (x) +2\Re\left(\sum \frac{\partial \Phi}{\partial
  z_k}(x)(z_k-z_k(x)) \right) &\geq&  \Phi (z) 
-\frac{1}{2} M_2  (2n)^2 \vert z-x\vert_{eucl}^2\\
&\geq&  \Phi (z) 
-\frac{1}{2} M_2  (2n)^2 \frac{r^2}{\mu ^2}.
\end{eqnarray*}

Then, for all  holomorphic
section $f$ of $L$ on $B_\om (p,r)$ and all $x$ in $B_\om (p,r/2)$,

\begin{eqnarray*}
  \lefteqn{\n f(x)\n_h^2= \vert f(x)\vert^2 e^{-\Phi (x)}}\\
&\leq& C
\int_{B_\om (x,\frac{r}{2})} 
\left\vert f e ^ { -\sum\limits_{1\leq\alpha\leq N} 
\overline{\sigma_\alpha (x)}(\sigma_\alpha (z)-\sigma_\alpha (x))-
\sum\limits_{1\leq k\leq n}  \frac{\partial \Phi_{def}}{\partial
    z_k}(x)(z_k-z_k(x))} \right\vert^2 e^{-\Phi (x)} dV_\om \\
&\leq& C  \int_{B_\om (x,r/2)} \n f(x)\n_h^2 dV_\om .
\end{eqnarray*}

\subsection{Gradient estimates}
Let $f$ be a holomorphic section of $L$ on $B_\om (p,r)$
and $x$ in $B_\om (p,r/4)$.
Denote by $s$ any coordinate among the $z_k-z_k(x)$ centered at $x$.
Then

\begin{eqnarray*}
 \lefteqn{ \frac{\partial}{\partial s} (\n f(x)\n_h^2)_{\vert s=0}}\\
 &=&
\overline{f(x)} e ^{-\Phi (x)}\left( \frac{\partial f}{\partial
    s}_{\vert s=0}-\frac{\partial \Phi}{\partial s}_{\vert s=0}
\right) \\
&=&
\overline{f(x)} e ^{-\Phi (x)} \frac{\partial }{\partial s}\left(f(s)e ^
{-\sum\limits_{1\leq\alpha\leq N} 
\overline{\sigma_\alpha (x)}(\sigma_\alpha (s)-\sigma_\alpha (x))
-\frac{\partial \Phi_{def}}{\partial s}(0)s}\right)_{\vert s=0}\\
&=& \overline{f(x)} e ^{-\frac{\Phi (x)}{2}} 
\int\limits_{\vert s\vert =\frac{r}{4\lambda}} 
f(s)e ^{-\sum\limits_{1\leq\alpha\leq N} 
\overline{\sigma_\alpha (x)}(\sigma_\alpha (s)-\sigma_\alpha (x))
-\frac{\partial \Phi_{def}}{\partial s}(0)s}
e ^{-\frac{\Phi (x)}{2}}
\frac{ds}{s}.
\end{eqnarray*}

Now, note that thanks to the regularity assumption
$$\left\vert f(s) e ^{-\sum\limits_{1\leq\alpha\leq N} 
\overline{\sigma_\alpha (x)}(\sigma_\alpha (z)-\sigma_\alpha (x))
-\frac{\partial \Phi_{def}}{\partial s}(0)s}\right\vert ^2 e ^
{-\Phi (x)}\leq C \n f(s)\n _h^2.$$

By the mean value estimate, we get 
$$\left\vert \frac{\partial}{\partial s} (\n f(x)\n_h^2)_{\vert s=0}
\right\vert   \leq C \int_{B_\om (p,r)} \n f\n_h^2 dV_\om.$$

\subsection{Uniform separation}
Choose an interpolating set $\Lambda$ for $\fc_{\om, h}$. Assume it is
not uniformly $\om$-separated. Then, there exists a subset $\Lambda '$ (still
interpolating) non-uniformly separated but such that any intersection
of three balls among the $B_\om (p_i,r)$ is empty. Hence, by the mean
value estimate, the evaluation map
$$\begin{array}{cccc}
  ev_{\Lambda '} :& \fc_{\om, h} &\to & l^2_{\Lambda ',h}\\
                  & f&           \mapsto & f_{\vert \Lambda '}
\end{array}
$$ is well defined and continuous.
By the open mapping theorem, $\Lambda '$ is stably interpolating, i.e.
$$\exists K, \forall a\in l^2_{\Lambda ',h}, \exists f\in \fc_{\om, h}
/ f_{\vert \Lambda '} =a \text{ and } \int_X \n f\n_h^2 dV_\om
\leq K \sum_{p\in\Lambda '}|a(p)|_h^2 <+\infty.$$ 
 
For all $p_0$ in $\Lambda '$ denote by $f^{p_0}$ any stable solution
for $a\not=0$ at ${p_0}$ and $a=0$ on $\Lambda '-\{p_0\}$.
 By the gradient estimate, we infer if $p\in\Lambda'\cap B_\om (p_0, r/4)$,
$$ \n f^{p_0}(p_0)\n_h^2 =\left\vert \n f^{p_0}(p_0)\n_h^2 -\n
  f^{p_0}(p)\n_h^2 \right\vert \leq 
\vert p-p_0 \vert  \int_X \n f\n_h^2 dV_\om \leq \frac{K}{\mu} d_\om
(p,p_0) \n f^{p_0}(p_0)\n_h^2 $$ which contradicts the assumption that
$\Lambda'$ is not uniformly separated.

\section{Interpolation on non-positively curved manifolds}
\label{sec:non-positive}

We will prove theorem~\ref{theo1} in this section.
We first have to describe the admissible metrics.
\begin{definition}
An hermitian line bundle $(L,h)\to U$ on an open set $U$ of $\cb^ n$
 has an $M_2$-almost analytic frame
$e$ on $U$ if the weight function $\Phi:= -\log h(e)$ 
can be written as the sum 
$$\sum_{1\leq \alpha\leq N} \vert \sigma_\alpha\vert ^2 + \Phi_{def}$$
where the $\sigma_\alpha$ are holomorphic functions on $U$ and 
all second order real derivatives of the deformation weight $\Phi_{def}$
are bounded by $-M_2$ and  $M_2$ on $U$.
\end{definition}
Note that nothing is now assumed on the oscillation of the 
holomorphic functions $\sigma_\alpha$.

\begin{definition}\label{defi:almost analytic}
An hermitian line bundle $(L,h)\to X$ is said to be
almost analytic around $\Lambda$  if there exists positive numbers $r_0$
and $M_2$ such that for all $p\in\Lambda$, there exists a coordinates chart 
$\Psi_p : (U,0)\subset \cb^n\to (V,p)\subset X$ centered at $p$ such
that 
\begin{itemize}
\item[(I)] $\Psi _p(U)\supset B_\om (p,r_0)$

\item[(II)] $(L,h)$ has a $M_2$-almost analytic frame $e$ on
  $\Psi_p^{-1}(B_\om (p,r_0))$. 

\item[(III)]  The differential $d\Psi_p :(TU,\Psi_p^\star \om (p)) \to
  (TV,\om)$ satisfies $\mu \leq \n d\Psi_p \n$ \\ on
  $\Psi_p^{-1}(B_\om (p,r_0))$. 
\end{itemize}
\end{definition}

\subsection{Local choice}
We now make the assumptions of theorem~\ref{theo1}. 
Denote by $2\delta_0 :=\inf _{{p,q\in\Lambda}\atop { p\not =q}} 
(d_\om(p ,q),r_0)$.
Let $a\in l^2_{\Lambda, h}$ and $p\in\Lambda$. 
On the ball $B_\om (p,\delta_0 )$, consider the holomorphic function
$$f_p := a(p) \exp \left(\sum\limits _{1\leq k\leq n}
 \frac{\partial \Phi_{def}}{\partial
    z_k}(p) (z_k - z_k (p))
+\sum\limits _{1\leq \alpha\leq N} \overline{\sigma_\alpha(p)}
  (\sigma_\alpha - \sigma_\alpha (p)) \right).$$
It is in fact a normal frame at $p$ for $(L,h)$.
The following estimates relies on the simple inequalities for $z$ in
$B_\om (p,\delta_0)$ which may be regarded as higher dimensional
analogues of Riesz formula.
\begin{eqnarray*}
  -\vert \sigma_\alpha (z)\vert ^2 +2\Re \left(
\overline{\sigma_\alpha (p)}(\sigma_\alpha (z)-\sigma_\alpha (p))\right)
& = &
-\vert \sigma_\alpha (p)\vert ^2 
- \vert \sigma_\alpha (z)-\sigma_\alpha (p)\vert ^2 \\
&\leq & -\vert \sigma_\alpha (p)\vert ^2 .
\end{eqnarray*}
and
\begin{eqnarray*}
  -\Phi_{def} (z) +2\Re\left(\sum\limits_{1\leq k\leq n}
 \frac{\partial \Phi_{def}}{\partial
  z_k}(p)(z_k-z_k(p)) \right)\!\!\!\! &\leq&\!\!\!\! -\Phi_{def} (p) 
+\frac{1}{2} M_2  (2n)^2 \vert z-p\vert_{eucl}^2\\
\!\!\!\! &\leq&\!\!\!\! -\Phi_{def} (p) 
+\frac{1}{2} M_2  (2n)^2 \frac{\delta_0^2}{\mu ^2}
\end{eqnarray*}

Then, for all $z\in B_\om (p,\delta_0 )$, 
\begin{eqnarray*}\label{eq}
\n f_p(z)\n_h^2=|f_p(z)|^2 e ^ {-\Phi (z)} &\leq&
C \n a(p)\n_h^2
\end{eqnarray*}
where $C$ is a constant independent of $p\in\Lambda$.

\subsection{Gluing}
 Let $\chi$ be a cut-off function defined on $\rb ^+$ decreasing,
 equal to 1  on  $[0,1/4]$, and to 0 on $[1,+\infty[$.
Then, $$F:= \sum _{p\in\Lambda} 
f_p \chi \left( \frac{d_\om (p,\cdot )^2}{\delta_0^2} \right)$$ 
is a smooth solution of the interpolation problem, holomorphic around 
$\Lambda$, with the growth condition 
\begin{eqnarray*}
\int_{X} \n F\n_h^2 dV_\om 
&\leq & \sum_{p\in\Lambda }\int_{B(p,\delta_0 )} C\n a(p)\n_h^2 dV_\om\\
&\leq&
C \sum_{p\in\Lambda }\n a(p)\n_h^2 <+\infty .
\end{eqnarray*}
Here, we used Rauch theorem to infer that 
the balls $B(p,\delta_0 )$ have volume bounded by a constant 
independent of $p$.

\subsection{Finding holomorphic sections}
We will apply the $L^2$-estimates existence theorem 
with a singular weight (see for example~\cite{de} theorem 5.1) 
in order to keep the right values on $\Lambda$.
Consider the auxiliary weight
$$v(z):=n \sum_{q\in\Lambda} \left(1- \frac{d_\om (q,z)^2}{\rho^2} + 
\log \frac{d_\om (q,z)^2}{\rho^2}\right) 1_{B_\om (q,\rho )}(z)$$

Remark first that,
\begin{eqnarray*}
\int_{X} \n\overline{\partial}F\n_h^2 e^{-v} dV_\om &= & 
\int_{X} \sum _{p\in\Lambda}  \left | \overline{\partial}\chi
\left(\frac{d_\om(p,z)^2}{\delta_0^2}\right) \right |^2 
\n f_p\n_h ^2e^{-v} dV_\om 
\end{eqnarray*} 
The term $|\overline{\partial}\chi
\left(\frac{d_\om(p,z)^2}{\delta_0^2}\right) |^2$
depends on the differential of the exponential map.
Hence it is bounded by a constant  independent of the point $p$ by
Rauch theorem. Now,
\begin{eqnarray*}
\int_{X} \n\overline{\partial}F\n_h^2 e^{-v} dV_\om 
&\leq & 
C \sum_{p\in\Lambda}
\n a(p) \n_h^2 \int_{B_\om (p,\delta _0)-B_\om (p,\delta _0 /2)}
e ^{-v}dV_\om .
\end{eqnarray*} 

We have to estimate $v$ around $p\in\Lambda$. Choose $q$ in $\Lambda$.
For $z\in (B_\om(p,\delta _0)-B_\om(p,\delta _0 /2))\cap B_\om(q,\rho)$, 
$d_\om (q,z)\geq d_\om (p,q) -d_\om (z,p) \geq 2\delta_0 -\delta_0/2$
 if $q\not =p$, and $d_\om (q,z) \geq \delta _0 /2$ if $q=p$. Hence,
$$1- \frac{d_\om (q,z)^2}{\rho^2} + \log \frac{d_\om (
  q,z)^2}{\rho^2}  \geq \log
\frac{d_\om(q,z)^2}{\rho^2} \geq \log \frac{\delta ^2_0}{ 4\rho ^2}$$
Recalling that the number of $q$ in $\Lambda$ such that $z\in B_\om(q,\rho)$ is
uniformly bounded because $\Lambda$ is assumed to be uniformly
separated and the volume of $B_\om(z,\rho)$ is bounded from above by
the volume of a ball radius $\rho$ in the hyperbolic space of
dimension $n$ and of sectional curvature $-k^2$, one infers that
 $$\int_{X} \n\overline{\partial}F\n_h^2 e^{-v} dx <+\infty.$$

Now, the computation of the curvature is simple because the 
derivatives of $$\left(1- \frac{d_\om (p,z)^2}{\rho^2} + \log
  \frac{d_\om (p,z)^2}{\rho^2}\right)$$ vanishes along $\partial
B_\om (p,a)$.
\begin{eqnarray*}
\hess v &=& n\sum _{p\in\Lambda} \left(\frac{-\hess 
  d_\om (p,z )^2 }{ \rho^ 2} +\hess \log {d_\om
  (p,z)^2}\right) 1_{B(p,a)}(z). 
\end{eqnarray*}
By the comparison theorem with the model $M=(\cb ,\vert~~\vert^2)$ and
$f=\log$ we get that $\log {d_\om (z,p)^2}$ is a
pluri-sub-harmonic function. 
On the other hand, using the hyperbolic
space of curvature $-k^2$ as model, with $f=x^2$, we get
(see~\cite{g-w} page 35)
$$\hess d_\om (p,z )^2 \leq (1+k d_\om (p,z ) \coth k
d_\om (p,z ) )\om \leq (1+k\rho\coth k\rho)\om$$
on $B_\om (p,\rho)$. 
Hence,
\begin{eqnarray*}
\lefteqn{ic_h(L)+ \hess v + ricci(\om )}\\
&\geq&ic_h(L)+ ricci(\om ) - n\frac{\card B_\om (z,\rho)\cap \Lambda}{\rho^2}
(1+k\rho\coth k\rho) \om\\
&\geq&\ep\om .
\end{eqnarray*}

Solving the equation $\overline{\partial}G=\overline{\partial}F$ on
$X$ endowed with the complete metric $\om$ for the
line bundle $L$  endowed with the 
metric $h e ^{-v}$, we get a smooth section $G$ of $L$ on $X$
such that $F-G$ is holomorphic and $\int_{X} \n G\n_h^2 e^{-v} dV_\om$ 
is finite. 
Now, because $v$ is non-positive,
$\int_{X} \n G\n_h^2 dV_\om$ is also finite. 
Furthermore, the estimate of $v$ around $p\in\Lambda$ leads to
$$\int_{B_\om(p,\delta _0)}\frac{\n G\n_h^2}{\n z-p\n^{2n}}\leq C
\int_{X} \n G \n_h^2 e^{-v} dx <+\infty  ,$$ so that $G$ vanishes
on $\Lambda$.

\subsection{Comments on the factor $\rho ^2$}

The first remark is an homogeneity property.
The quantity $$\displaystyle\frac{\card B_\omega (z,\rho)\cap
  \Lambda}{\rho^2}\omega $$
which appears in theorem~\ref{theo1} applied for flat manifolds is equal to 
the quantity $\displaystyle\frac{\card B_{t\omega}(z,\sqrt{t}\rho)\cap
  \Lambda}{(\sqrt{t}\rho)^2}t\omega$ for all positive $t$.

The second remark is a restriction property, which, together with the
criterion for interpolation in $\cb$ yields a reason for the factor 
$\rho^ 2$.
\begin{prop}
Consider $\cb^n$ with its usual flat metric $\om$ and $(L,h)\to \cb^n$
an hermitian line bundle with $M_2$-regular frame on a neighbourhood of 
$\cb^{n-1}$ in $\cb^n$.

 If $\Lambda$ is
interpolating for $\fc _{\om,h}$ then $\Lambda\cap \cb^{n-1}$ is
interpolating for $\fc _{\om _{| \cb ^{n-1}}, h_{ | \cb ^{n-1}}}$. 
\end{prop}

\textbf{Proof}
Let $a\in l^2_{\Lambda \cap \cb^{n-1},h}$.
There exists a function $f$ in $\fc_{\om,h}$ achieving $f(p)=a(p)$ for
$p\in\Lambda\cap \cb^{n-1}$ and $f(p)=0  $ for
$p\in\Lambda - \cb^{n-1}$. To check the growth condition, 
\begin{eqnarray*}
\lefteqn{|f(z',0)|^2 e ^{-\Phi (z',0)}}\\ &\leq &\!\!\!\!
\frac{1}{\pi\delta_0^2}\int\limits_{z_n\in D(0,\delta_0)} \left| f(z',z_n)e ^
{-\sum\limits_{1\leq\alpha\leq N} 
\overline{\sigma_\alpha (z',0)}(\sigma_\alpha (z)-\sigma_\alpha (z',0))
-\frac{\partial\Phi_{def}}{\partial z_n}(z',0)z_n }\right| ^2 e
^{-\Phi(z',0)} \\
& \leq &\!\!\!\!
C e ^ {1/2 M_2 2^2\delta_0^2 + NM_2}
\int\limits_{z_n\in D(0,\delta_0)} |f(z',z_n)|^2 e ^{-\Phi(z',z_n)} dx_n dy_n.
\end{eqnarray*}

We can now apply Fubini theorem because the measure on $\cb^n$ is a
product measure. We get

$\int_{\cb^{n-1}} |f(z',0)|^2 e ^{-\Phi(z',0)} \leq
C\int_{\cb ^n} |f(z)|^2 e ^{-\Phi(z)} <+\infty.$ 

\section{Interpolation on negatively curved manifolds}
\label{sec:negative}
We will now prove theorem~\ref{theo2}. 

The ball $\bb (\kappa )$ of radius $\kappa$ in $\rb^{2n}$ endowed with
the metric $\displaystyle \frac{4\vert dz\vert^2}{(1-\frac{\vert
    z\vert^2}{\kappa ^2})^2}$ has constant sectional curvature
$-1/\kappa ^2$. The distance between $0$ and $z$ is $2\kappa \text{arctanh}
\frac{\vert z\vert}{\kappa}$. Hence,
$$\hess \log \tanh^2 \frac{d_{\bb (\kappa )} (0,\cdot )}{2\kappa} = \hess
\log  \vert z\vert ^2 \geq 0.$$
Hence, by comparison theorem, for every $p\in X$,
$$\hess \log \tanh^2\frac{d_\om (p,\cdot )}{2\kappa} \geq 0.$$

Moreover, $v_p :=n\log \tanh^2\frac{d_{X} (p,\cdot )}{2\kappa}$ is negative and
  has a pole of order $2n$ at $p$. Consider $v:=\sum_{p\in\Lambda}
  v_p$. By assumption on the density, 
$$v(z)\geq -\sup_X D_{\Lambda ,\kappa} + \sum_{p\in\Lambda\atop d_\om (z,p)\leq
  1} v_p(z) $$
Just as we estimated the auxiliary weight in the previous proof, we
can estimate $v$ on  $B_\om(p,\delta _0)-B_\om(p,\delta _0 /2)$.
The rest of the proof follows the same lines as the proof of 
theorem~\ref{theo1}.

The function $v$ may be compared to a pluri-complex Green function
having poles on $\Lambda$.

\end{document}